\documentclass[12pt]{article}
\usepackage{mathrsfs}
\usepackage{amsmath,amsfonts,amssymb,rotating,amsthm}
\usepackage{hyperref}
\usepackage{array}

\def\qed{\nopagebreak\hfill{\rule{4pt}{7pt}}}
\def\proof{\noindent {\it{Proof.} \hskip 2pt}}

\newtheorem{theo}{Theorem}[section]

\newtheorem{prop}[theo]{Proposition}

\newtheorem{conj}[theo]{Conjecture}

\theoremstyle{remark}

\numberwithin{equation}{section}

\allowdisplaybreaks 

\newdimen\Squaresize \Squaresize=11pt
\newdimen\Thickness \Thickness=0.7pt
\def\Square#1{\hbox{\vrule width \Thickness
   \vbox to \Squaresize{\hrule height \Thickness\vss
    \hbox to \Squaresize{\hss#1\hss}
   \vss\hrule height\Thickness}
\unskip\vrule width \Thickness} \kern-\Thickness}

\def\Vsquare#1{\vbox{\Square{$#1$}}\kern-\Thickness}

\def\moins{\raise 1pt\hbox{{$\scriptstyle -$}}}

\begin{document}

\begin{center}
{\large \bf  Finite Differences of the Logarithm of

the Partition Function}
\end{center}

\begin{center}
William Y.C. Chen$^{1}$, Larry X.W. Wang$^{2}$ and Gary Y.B. Xie$^{3}$ \\[8pt]
Center for Combinatorics, LPMC-TJKLC\\
Nankai University\\
 Tianjin 300071, P. R. China\\[6pt]
Email: $^{1}${\tt chen@nankai.edu.cn}, $^{2}${\tt
wsw82@nankai.edu.cn, $^{3}$xieyibiao@mail.nankai.edu.cn}
\end{center}

\vspace{0.3cm} \noindent{\bf Abstract.} Let $p(n)$ denote the partition function.
DeSalvo and Pak proved that  $\frac{p(n-1)}{p(n)}\left(1+\frac{1}{n}\right)> \frac{p(n)}{p(n+1)}$ for $n\geq 2$, as conjectured by Chen.
Moreover, they conjectured that a sharper inequality $\frac{p(n-1)}{p(n)}\left( 1+\frac{\pi}{\sqrt{24}n^{3/2}}\right) > \frac{p(n)}{p(n+1)}$ holds for $n\geq 45$. In this paper, we prove the conjecture of Desalvo and Pak  by giving an upper bound for $-\Delta^{2} \log p(n-1)$, where $\Delta$ is the difference operator with respect to $n$. We also show that for given $r\geq 1$ and sufficiently large $n$, $(-1)^{r-1}\Delta^{r} \log p(n)>0$. This is analogous to the
 positivity of finite differences of the partition function.
 It was conjectured by Good and proved by Gupta  that for given $r\geq 1$, $\Delta^{r} p(n)>0$ for sufficiently large $n$.

\noindent {\bf Keywords:} partition function, log-concavity, finite difference, the
Lambert $W$ function, the Hardy-Ramanujan-Rademacher formula

\noindent {\bf AMS Classification:} 05A20, 11B68

\section{Introduction}

A partition of positive integer $n$ is a nonincreasing sequence of positive integers $\lambda_{1},\lambda_{2},\ldots,\lambda_{r}$ such that $\sum_{i=1}^{r} \lambda_{i}=n$. Let $p(n)$ denote the number of partitions of $n$. In particular, we set $p(0)=1$. The Hardy-Ramanujan-Rademacher formula for $p(n)$ states that
\begin{equation*}
\begin{split}
  p(n)=&\frac{\sqrt{12}}{24n-1}\sum^{N}_{k=1}A_{k}(n)\sqrt{k}
  \left[\left(1-\frac{k}{\mu(n)}\right)e^{\mu(n)/k}
  +\left(1+\frac{k}{\mu(n)}\right)e^{-\mu(n)/k}\right]\\
  &\quad +R_{2}(n,N),
\end{split}
\end{equation*}
where $A_{k}(n)$ is an arithmetic function, $R_{2}(n,N)$ is the remainder term and
\begin{equation}\label{mumuref1}
 \mu(n)=\frac{\pi}{6}\sqrt{24n-1},
\end{equation}
 see, for example, Hardy and Ramanujan \cite{hrsf}, Rademacher \cite{rade}. Note that $A_{1}(n)=1$ and $A_{2}(n)=(-1)^{n}$ for $n\geq 1$. Lehmer \cite{lehmm,lehm} gave the following error bound
\begin{equation*}
  |R_{2}(n,N)|<\frac{\pi^{2}N^{-2/3}}{\sqrt{3}}\left[\left(\frac{N}{\mu(n)}\right)^{3}
  \sinh \frac{\mu(n)}{N}+\frac{1}{6}-\left(\frac{N}{\mu(n)}\right)^{2}\right],
\end{equation*}
which is valid for all positive integers $n$ and $N$.

Employing Rademacher's convergent series and Lehmer's error bound, DeSalvo and Pak \cite{desp} proved the following inequality conjectured by Chen \cite{chen3}.

\begin{theo}
For $n\geq 2$, we have
\begin{equation}
\frac{p(n-1)}{p(n)}\left( 1+\frac{1}{n}\right)> \frac{p(n)}{p(n+1)}.
\end{equation}
\end{theo}

The above relation has been improved by DeSalvo and Pak \cite{desp}.

\begin{theo}
For $n \geq 7$, we have
\begin{equation}\label{dpconj1}
\frac{p(n-1)}{p(n)}\left( 1+\frac{240}{(24n)^{3/2}}\right)> \frac{p(n)}{p(n+1)}.
\end{equation}
\end{theo}

They also proposed the following conjecture.

\begin{conj}\label{dp}
For $n\geq 45$, we have
\begin{equation}\label{dpc1}
\frac{p(n-1)}{p(n)}\left( 1+\frac{\pi}{\sqrt{24}n^{3/2}}\right)> \frac{p(n)}{p(n+1)}.
\end{equation}
\end{conj}

It should be mentioned that by using Lehmer's error bound for the remainder term of $p(n)$,  Bessenrodt and Ono \cite{bo} proved the following inequality.

\begin{theo}
For any integers $a,b$ satisfying $a,b>1$ and $a+b>9$, we have
\[
p(a)p(b)>p(a+b).
\]
\end{theo}

In this paper,  we shall prove Conjecture \ref{dp}  by giving an upper bound for $-\Delta^2 \log p(n-1)$ for $n\geq 5000$. Moreover, for any given $r$, we give an upper bound for $(-1)^{r-1}\Delta ^{r} \log p(n)$.

In 1977, Good \cite{goodp} conjectured that $\Delta^{r}p(n)$ alternates in sign up to a certain value $n=n(r)$, and then it stays positive. Using the Hardy-Rademacher series \cite{rade1} for $p(n)$, Gupta \cite{gupta} proved that for any given $r$, $\Delta^{r} p(n)>0$ for sufficiently large $n$. In 1988, Odlyzko \cite{amo} proved the conjecture of Good and obtained the following asymptotic formula for $n(r)$:
\begin{equation*}
  n(r)\sim \frac{6}{\pi^{2}}r^{2}\log^{2} r \quad \textrm{as}~ r \rightarrow \infty.
\end{equation*}
Knessl and Keller \cite{kka1,kka} obtained an approximation $n(r)'$ for $n(r)$ for which  $|n(r)'-n(r)|\leq 2$ up to $r=75$.  Almkvist \cite{gadp1, gadp} proved that $n(r)$ satisfies certain equations.

By using the bounds of the modified Bessel function of the first kind, we shall prove that for any given $r\geq 1$, there exists a positive integer $n(r)$ such that $(-1)^{r-1}\Delta ^{r} \log p(n)>0$ for $n\geq n(r)$.

\section{Proof of Conjecture \ref{dp}}

In this section, we give a proof of Conjecture \ref{dp} by using an inequality of DeSalvo and Pak \cite{desp}.
Let \[p_{2}(n)=2\log p(n)-\log p(n-1)-\log p(n+1),\] DeSalvo and Pak have shown that  for $n \geq 50$,
\begin{eqnarray}\label{p22}
  p_{2}(n) \!&<&\! \frac{24\pi}{(24(n-1)-1)^{3/2}}+\frac{288\pi(-3+\pi \sqrt{24(n-1)-1})}{(24(n-1)-1)^{3/2}(-6+\pi \sqrt{24(n-1)-1})^{2}} \nonumber\\
   & & \quad -\frac{864}{(24(n+1)-1)^{2}}+2e^{-\frac{\pi}{10}\sqrt{\frac{2n}{3}}}.
\end{eqnarray}
We shall give an estimate of the right hand side of (\ref{p22}), leading to a
proof of the conjecture.

\noindent {\emph{Proof of Conjecture} \ref{dp}.} The conjecture can be
restated as follows
\begin{equation}\label{cpp1}
 p_2(n)< \log \left(1+\frac{\pi}{\sqrt{24}n^{3/2}}\right),
\end{equation}
where $n\geq 45$. We proceed to give an estimate of each term of
the right hand side of (\ref{p22}).

We begin with the first term of the right hand side of \eqref{p22}.
We claim that for $n\geq 50$,
\begin{equation}\label{cpp3}
\frac{24\pi}{(24(n-1)-1)^{3/2}}<\frac{24\pi}{(24n)^{3/2}}-
\left(\frac{24\pi}{(24n)^{3/2}}\right)^{2}+\frac{3}{2n^{5/2}}.
\end{equation}
For $0<x\leq \frac{1}{48}$, it can be easily checked that
\begin{equation}\label{cpp2}
 \frac{1}{(1-x)^{3/2}}< 1+\frac{3}{2}x+\frac{3}{8}x^{3/2}.
\end{equation}
For $n\geq 50$, we have $\frac{25}{24n}\leq \frac{1}{48}$, and hence
we can apply (\ref{cpp2}) to deduce that
\begin{eqnarray}\label{hhh}
\frac{24\pi}{(24(n-1)-1)^{3/2}}&=&
  \frac{24\pi}{(24n)^{3/2}\left(1-\frac{25}{24n}\right)^{3/2}}\nonumber\\
  &<&\frac{24\pi}{(24n)^{3/2}}
     \left(1+\frac{3}{2}\cdot\frac{25}{24n}+\frac{3}{8}
     \left(\frac{25}{24n}\right)^{3/2}\right).
\end{eqnarray}
For $n\geq 50$, we have
\begin{eqnarray*}
\frac{3}{8}\left(\frac{25}{24n}\right)^{3/2}&<&
\frac{3}{8}\left(\frac{25}{24}\right)^{3/2}
\cdot \frac{1}{50^{1/2}n},\\[6pt]
\frac{24\pi}{(24n)^{3/2}}&<&\frac{24\pi}{(24)^{3/2}50^{1/2}}\cdot \frac{1}{n}.\\
\end{eqnarray*}
It follows that
\begin{eqnarray}\label{dpak1}
  &&\frac{24\pi}{(24n)^{3/2}}\left(\frac{3}{2}\cdot\frac{25}{24n}+
  \frac{3}{8}\left(\frac{25}{24n}\right)^{3/2}+\frac{24\pi}{(24n)^{3/2}}\right)
  \nonumber\\[6pt]
  &&\quad \leq \frac{24\pi}{(24n)^{3/2}}\cdot\frac{1}{n}\left(\frac{25}{16}+
  \frac{3}{8}\left(\frac{25}{24}\right)^{3/2}\cdot \frac{1}{50^{1/2}}+\frac{24\pi}{(24)^{3/2}50^{1/2}}\right)\nonumber\\
  &&\quad <\frac{3}{2n^{5/2}} .
\end{eqnarray}
Combining (\ref{hhh}) and (\ref{dpak1}), we obtain (\ref{cpp3}).

As for the second term of the right hand side of \eqref{p22}, it can be shown that for $n>50$,
\begin{equation}\label{cpp5}
\frac{288\pi(-3+\pi\sqrt{24(n-1)-1})}{(24(n-1)-1)^{3/2}(-6+\pi\sqrt{24(n-1)-1})^{2}}
<\frac{1}{2n^{2}}+\frac{1}{n^{5/2}}.
\end{equation}
To this end, we need the following inequality for $\alpha \geq \frac{1}{2}$ and $0< x\leq c<1$,
\begin{equation}\label{gongs1}
\frac{1}{(1-x)^{\alpha}}\leq 1+\left(\frac{1}{1-c}\right)^{\alpha+1}\alpha x.
\end{equation}
Let
 \[f(x)=\frac{1}{(1-x)^{\alpha}}- 1-\left(\frac{1}{1-c}\right)^{\alpha+1}\alpha x.\]
For $\alpha \geq \frac{1}{2}$ and $0\leq x\leq c<1$, we see that
\begin{equation*}
f'(x)=\frac{\alpha}{(1-x)^{\alpha+1}}-\left(\frac{1}{1-c}\right)^{\alpha+1}\alpha \leq 0.
\end{equation*}
Since $f(0)=0$, we obtain that $f(x)\leq 0$ under the above assumption.
This yields that $f(x)<0$ for $0< x\leq c<1$ and $\alpha \geq \frac{1}{2}$,
and hence (\ref{gongs1}) is proved.

The left hand side of \eqref{cpp5} can be rewritten as
\begin{equation*}
  \frac{288\pi\cdot \frac{\pi}{2} \sqrt{24n-25}}{(24n-25)^{3/2}(-6+\pi\sqrt{24n-25})^{2}}+
  \frac{288\pi(-3+\frac{\pi}{2}\sqrt{24n-25})}{(24n-25)^{3/2}(-6+\pi\sqrt{24n-25})^{2}},
\end{equation*}
which can be simplified to
\begin{equation}\label{dpak3}
  \frac{1}{4n^{2}\left(1-\frac{25}{24n}\right)^{2}
  \left(1-\frac{6}{\pi\sqrt{24n-25}}\right)^{2}}+
  \frac{1}{4n^{2}\left(1-\frac{25}{24n}\right)^{2}\left(1-\frac{6}{\pi\sqrt{24n-25}}\right)}.
\end{equation}
Setting $x =\frac{25}{24n}$, $\alpha=2$ and $c=\frac{1}{48}$,
for $n\geq 50$, we have $0<x<c<1$ and $\alpha\geq \frac{1}{2}$. By (\ref{gongs1}), we find that for $n\geq 50$,
\begin{equation}\label{dpak4}
  \frac{1}{\left(1-\frac{25}{24n}\right)^{2}}\leq  1+\left(\frac{48}{47}\right)^{3}\frac{25}{12n}.
\end{equation}
Setting $x=\frac{6}{\pi\sqrt{24n-25}}$, $\alpha=2$ and $c=\frac{1}{15}$, for
$n\geq 50$, we also have $0<x<c<1$ and $\alpha \geq \frac{1}{2}$. Again, using
(\ref{gongs1}), we see that for $n\geq 50$,
\begin{equation}\label{dpak5}
  \frac{1}{\left(1-\frac{6}{\pi\sqrt{24n-25}}\right)^{2}}<
  1+\left(\frac{15}{14}\right)^{3}\frac{6}{\pi \sqrt{24n-25}}<1+\frac{24}{\pi \sqrt{24n-25}}.
\end{equation}
Combining \eqref{dpak4} and \eqref{dpak5}, we deduce that for $n\geq 50$,
\begin{eqnarray}\label{dpak6}
  &&\frac{1}{4n^{2}\left(1-\frac{25}{24n}\right)^{2}
  \left(1-\frac{6}{\pi\sqrt{24n-25}}\right)^{2}}\nonumber\\ &&\quad \leq\frac{1}{4n^{2}}\left(1+\left(\frac{48}{47}\right)^{3}
  \frac{25}{12n}\right)\left(1+\frac{24}{\pi \sqrt{24n-25}}\right).
\end{eqnarray}
Setting $x =\frac{25}{24n}$, $\alpha=\frac{1}{2}$ and $c=\frac{1}{48}$,
for $n\geq 50$, we have $0<x<c<1$ and $\alpha\geq \frac{1}{2}$. Applying (\ref{gongs1}), for $n\geq 50$, we get
\begin{eqnarray}\label{dpak7}
\frac{24}{\pi \sqrt{24n-25}}&=&\frac{24}{\pi (24n)^{1/2}}\frac{1}{\left(1-\frac{25}{24n}\right)^{1/2}}\nonumber\\
&<&\frac{24}{\pi (24n)^{1/2}}\left(1+\left(\frac{48}{47}\right)^{3/2}\frac{25}{48n}\right).
\end{eqnarray}
Combining \eqref{dpak6} and \eqref{dpak7}, we find that  for $n\geq 50$,
\begin{align}
&\frac{1}{4n^{2}\left(1-\frac{25}{24n}\right)^{2}
\left(1-\frac{6}{\pi\sqrt{24n-25}}\right)^{2}}\notag\\
& \quad \leq \frac{1}{4n^{2}}\left(1+\left(\frac{48}{47}\right)^{3}\frac{25}{12n}\right)
\left(1+\frac{24}{\pi (24n)^{1/2}}\left(1+\left(\frac{48}{47}\right)^{3/2}\frac{25}{48n}\right)\right).
\label{dpak8}
\end{align}
The right hand side of \eqref{dpak8} can be expanded as follows
\begin{align}
&\frac{1}{4n^{2}}+\frac{\sqrt{6}}{2\pi n^{5/2}}+\frac{25}{48n^{3}}\left(\frac{48}{47}\right)^{3}
+\frac{25\sqrt{6}}{96\pi n^{7/2}}\left(\frac{48}{47}\right)^{3/2} \notag \\
&\quad +\frac{25\sqrt{6}}{24\pi n^{7/2}}\left(\frac{48}{47}\right)^{3}
+\frac{25^{2}\sqrt{24}}{48^{2}\pi n^{9/2}}\left(\frac{48}{47}\right)^{9/2}.\label{gcpp1}
\end{align}
Clearly, for $\alpha>\frac{5}{2}$ and $n\geq 50$,
\begin{equation*}  
\frac{1}{n^{\alpha}}\leq \frac{1}{50^{\alpha-5/2}n^{5/2}},
\end{equation*}
which implies that for $n\geq 50$,
\begin{equation}\label{examgp1}
\frac{1}{n^{3}}\leq \frac{1}{50^{1/2}n^{5/2}},
\end{equation}
\begin{equation}\label{examgp2}
\frac{1}{n^{7/2}}\leq \frac{1}{50n^{5/2}},
\end{equation}
\begin{equation}\label{examgp3}
\frac{1}{n^{9/2}}\leq \frac{1}{50^{2}n^{5/2}}.
\end{equation}
Applying \eqref{examgp1}, \eqref{examgp2} and \eqref{examgp3} to the last four terms of \eqref{gcpp1}, we obtain that for $n\geq 50$,
\begin{equation}\label{22s1}
\frac{1}{4n^{2}\left(1-\frac{25}{24n}\right)^{2}\left(1-\frac{6}{\pi\sqrt{24n-25}}\right)^{2}} <\frac{1}{4n^{2}}+\frac{1}{2n^{5/2}}.
\end{equation}
For the second term of \eqref{dpak3}.
Setting $x=\frac{6}{\pi\sqrt{24n-25}}$, $\alpha=1$ and $c=\frac{1}{15}$, for
$n\geq 50$, we  have $0<x<c<1$ and $\alpha \geq \frac{1}{2}$. By
(\ref{gongs1}), we see that for $n\geq 50$,
\begin{equation}\label{dpaks5}
  \frac{1}{1-\frac{6}{\pi\sqrt{24n-25}}}<1+\left(\frac{15}{14}\right)^{2}\frac{6}{\pi \sqrt{24n-25}}<1+\frac{12}{\pi \sqrt{24n-25}}.
\end{equation}
Using \eqref{dpaks5} and the same argument as in the derivation of
\eqref{22s1}, it can be shown that for $n\geq 50$,
\begin{equation}\label{22s2}
\frac{1}{4n^{2}\left(1-\frac{25}{24n}\right)^{2}\left(1-\frac{6}{\pi\sqrt{24n-25}}\right)}
<\frac{1}{4n^{2}}+\frac{1}{2n^{5/2}}.
\end{equation}
In view of (\ref{22s1}) and (\ref{22s2}), we arrive at \eqref{cpp5}.

To estimate the third term of the right hand side of \eqref{p22}, we aim
 to show that for $n\geq 50$,
\begin{equation}\label{cpp4}
-\frac{864}{(24(n+1)-1)^{2}}<\frac{1}{2n^{5/2}}-\frac{3}{2n^{2}}.
\end{equation}
It's easily verified that for $\alpha\geq 1/2$ and $0\leq x\leq1$,
\begin{equation}\label{posalpha1}
1\geq \frac{1}{(1+x)^{\alpha}}\geq 1-\alpha x.
\end{equation}
So for $n\geq 50$, we have
\[\frac{1}{\left(1+\frac{23}{24n}\right)^{2}}\geq 1-\frac{23}{12n}.\]
Consequently, for $n\geq 50$,
 \begin{equation*}
   -\frac{864}{(24(n+1)-1)^{2}}=-\frac{3}{2n^{2}\left(1+\frac{23}{24n}\right)^{2}}
   \leq\frac{23}{8n^{3}}-\frac{3}{2n^{2}}
   \leq \frac{1}{2n^{5/2}}-\frac{3}{2n^{2}}.
 \end{equation*}

Utilizing the above upper bounds \eqref{cpp3}, \eqref{cpp5} and \eqref{cpp4} for the three terms of the right hand side of (\ref{p22}), we conclude that for $n\geq 50$,
\begin{equation*}
  p_{2}(n)< \frac{24\pi}{(24n)^{3/2}}-\left(\frac{24\pi}{(24n)^{3/2}}\right)^{2}
  -\frac{1}{n^{2}}+\frac{3}{n^{5/2}}+2e^{-\frac{\pi}{10}\sqrt{\frac{2n}{3}}}.
\end{equation*}

Next we show that for $n\geq 5000$,
\begin{equation}\label{cwfcc2}
  p_{2}(n)<\frac{24\pi}{(24n)^{3/2}}-\left(\frac{24\pi}{(24n)^{3/2}}\right)^{2}.
\end{equation}
Clearly,  for $n\geq 100$, \[-\frac{1}{n^{2}}+\frac{3}{n^{5/2}}<-\frac{2}{3n^{2}}.\]
To prove that for $n\geq 5000$,
\begin{equation}\label{cwfc1}
  -\frac{2}{3n^{2}}+2e^{-\frac{\pi}{10}\sqrt{\frac{2n}{3}}}<0,
\end{equation}
let
\[ g(x)=-\frac{2}{3 x^{2}}+2e^{-\frac{\pi}{10}\sqrt{\frac{2x}{3}}}.\]
The equation $g(x)=0$ has two solutions
\begin{eqnarray*}
x_{1}&=&\frac{2400}{\pi^{2}}\left(W_{0}\left(-\frac{\pi\sqrt{2}}{40\cdot 3^{3/4}}\right)\right)^{2},\\
x_{2}&=& \frac{2400}{\pi^{2}}\left(W_{-1}\left(-\frac{\pi\sqrt{2}}{40\cdot 3^{3/4}}\right)\right)^{2},
\end{eqnarray*}
where $W_{0}(z)$ and $W_{-1}(z)$ are two branches of Lambert $W$ function $W(z)$, see Corless, Gonnet, Hare, Jeffrey and Knuth \cite{lwfk}.
More explicitly, we have $x_{1}\approx 0.64$, $x_{2}\approx 4996.47$.
It can be checked that $g(5000)<0$. Thus for $x\geq 5000$, \[g(x)<0.\]
This proves \eqref{cwfc1}. Hence (\ref{cwfcc2}) holds.

Using (\ref{cwfcc2}), it can be shown that (\ref{cpp1}) holds for $n\geq 5000$.
It is easily verified that for $x>0$,
\begin{equation}\label{cwfcc3}
x(1-x)< \log(1+x).
\end{equation}
Let
\begin{equation*}
h(x)=\log(1+x)-x+x^{2}.
\end{equation*}
For $x\geq 0$, we see that
\begin{equation*}
h'(x)=\frac{x+2x^{2}}{1+x}\geq 0.
\end{equation*}
Since $h(0)=0$, we have $h(x)>0$ for $x>0$.
Combining (\ref{cwfcc2}) and (\ref{cwfcc3}), we conclude that for $n \geq 5000$,
\begin{equation*}
p_2(n)< \log \left(1+\frac{\pi}{\sqrt{24}n^{3/2}}\right).
\end{equation*}
DeSalvo and Pak \cite{desp} have verified the above relation for $45\leq n\leq 8000$.
Thus (\ref{cpp1}) holds for $n\geq 45$ and the proof is completed. \qed

\section{An upper bound for $(-1)^{r-1} \Delta ^{r}\log p(n)$}

The conjecture of DeSalvo and Pak can be formulated as an
upper bound for $ 2\log p(n)-\log p(n-1)-\log p(n+1) $, namely, for $n\geq 45$,
\begin{equation}\label{betcpp1}
 -\Delta^{2}\log p(n-1) < \log \left(1+\frac{\pi}{\sqrt{24}n^{3/2}}\right),
\end{equation}
where $\Delta$ is the difference operator as given by $\Delta f(n)=f(n+1)-f(n)$.

 In this section, we  give an upper bound for $(-1)^{r-1} \Delta^{r} \log p(n)$.
 When $r=2$, this upper bound reduces to the above relation (\ref{betcpp1}).
 In the following theorem, we adopt the notation $(a)_k$ for the
 rising factorial, namely,    $(a)_0=1$ and $(a)_{k}=a(a+1)\cdots (a+k-1)$ for $k\geq 1$.

\begin{theo}\label{dprm1}
For each $r\geq 1$, there exists a positive integer $n(r)$ such that for $n\geq n(r)$,
\begin{equation*}
  (-1)^{r-1}\Delta
   ^{r} \log p(n)<\log \left(1+\frac{\sqrt{6}\pi}{6} \left(\frac{1}{2}\right)_{r-1}\frac{1}{(n+1)^{r-\frac{1}{2}}}\right).
\end{equation*}
\end{theo}

In the proof of the above theorem, we shall
use Hardy-Ramanujan-Rademacher series for $n\geq 1$,
\begin{equation}\label{hrs3}
  p(n)=2\pi \left(\frac{\pi}{12}\right)^{3/2}\sum_{k=1}^{\infty}A_{k}(n)k^{-5/2}L_{3/2}\left(\frac{\pi^2}{6k^2}\left(n-\frac{1}{24}\right)\right),
\end{equation}
and the following estimate for $A_{k}(n)$,
\begin{equation}\label{amr2}
    |A_{k}(n)|\leq 2k^{3/4},
\end{equation}
see Rademacher \cite{rade1}. In particular, we have $A_{1}(n)=1$ and $A_{2}(n)=(-1)^{n}$. The function $L_{\nu}(x)$ in (\ref{hrs3}) is defined by
\begin{equation}\label{deflnux1}
  L_{\nu}(x)=\sum_{m=0}^{\infty}\frac{x^{m}}{m!\Gamma(m+\nu+1)},
\end{equation}
where $\Gamma(m+\nu+1)$ is the Gamma function.

With the notation of $\mu(n)$ as in (\ref{mumuref1}), we have
 \[ \frac{\pi^2}{6}\left(n-\frac{1}{24}\right)=\frac{\mu^{2}(n)}{4},\]
 and so \eqref{hrs3} can be rewritten as
\begin{equation}\label{ccphrs3}
  p(n)=2\pi \left(\frac{\pi}{12}\right)^{3/2}\sum_{k=1}^{\infty}A_{k}(n)k^{-5/2}
  L_{3/2}\left(\frac{\mu^2(n)}{4k^2}\right),
\end{equation}
Denote the $k$th summand in \eqref{ccphrs3} by $f_{k}(n)$, namely,
\begin{equation}\label{askm1}
 f_k(n)=2\pi \left(\frac{\pi}{12}\right)^{3/2}A_{k}(n)k^{-5/2}
  L_{3/2}\left(\frac{\mu^2(n)}{4k^2}\right).
\end{equation}
Writing \eqref{ccphrs3} as
\begin{equation}\label{ccpheq3}
p(n)=f_{1}(n)\left(1+\frac{f_{2}(n)}{f_{1}(n)}\right)\left(1+\frac{\sum_{k\geq 3}^{\infty}f_{k}(n)}{f_{1}(n)+f_{2}(n)}\right).
\end{equation}
It is known that
\begin{equation*}
  L_{3/2}(x)=\frac{1}{\sqrt{\pi}}\cdot \frac{d}{dx}\left(\frac{\sinh 2\sqrt{x}}{\sqrt{x}}\right),
\end{equation*}
see Abramowitz and Stegun \cite{ahand} or Almkvist \cite{gadp1}.
Since $A_1(n)=1$,  $f_{1}(n)$ can be  expressed as
\begin{equation}\label{ccphrs4}
  f_{1}(n)=\frac{\sqrt{12}}{24n-1}
  \left[\left(1-\frac{1}{\mu(n)}\right)e^{\mu(n)}
  +\left(1+\frac{1}{\mu(n)}\right)e^{-\mu(n)}\right].
\end{equation}
Recall $A_2(n)=(-1)^n$, by \eqref{deflnux1} and \eqref{askm1}
we obtain that for $n\geq 1$,
\begin{equation*}
 f_{1}(n)-|f_{2}(n)|=2\pi \left(\frac{\pi}{12}\right)^{3/2}\sum_{m=0}^{\infty}
 \left(\frac{1}{4^m}-\frac{1}{2^{5/2}16^m}\right) \frac{\mu^{2m}(n)}{m!\Gamma(m+5/2)}.
\end{equation*}
Clearly, $\frac{1}{4^m}-\frac{1}{2^{5/2}16^m}>0$ for $m\geq 0$.
Hence for $n \geq 1$,
\begin{equation}\label{geq0}
f_{1}(n)-|f_{2}(n)|>0,
\end{equation}
which implies that for $n\geq 1$, $f_1(n)$ is positive and
\begin{equation*}
f_{1}(n)+f_{2}(n)>0.
\end{equation*}
It is also clear that, for $n\geq 1$, both of $\mu(n)-1$ and $1+\frac{\sum_{k\geq 3}^{\infty}f_{k}(n)}{f_{1}(n)+f_{2}(n)}$ are  positive.
Applying (\ref{ccphrs4}) to \eqref{ccpheq3}, we obtain that for $n\geq 1$
\begin{align*}
  \log p(n)=&\log \frac{\pi^{2}}{6\sqrt{3}}-3\log \mu(n) +\log (\mu(n)-1)+\mu(n)\notag\\
  &\quad+\log \left(1+\frac{\mu(n)+1}{\mu(n)-1}e^{-2\mu(n)}\right)+\log \left(1+\frac{f_{2}(n)}{f_{1}(n)}\right) \notag\\
  &\quad+\log \left(1+\frac{\sum_{k\geq 3}^{\infty}f_{k}(n)}{f_{1}(n)+f_{2}(n)}\right).
\end{align*}
Hence
\begin{equation}\label{loghrs2}
  (-1)^{r-1}\Delta^{r}\log p(n)=H_{r}+F_{1}+F_{2}+F_{3},
\end{equation}
where
\begin{eqnarray*}
  H_{r} &=& (-1)^{r-1}\Delta^{r}\left(-3\log \mu(n)
  +\log (\mu(n)-1)+\mu(n) \right), \\[6pt]
  F_{1} &=& (-1)^{r-1}\Delta^{r} \log \left(1+\frac{\mu(n)+1}{\mu(n)-1}e^{-2\mu(n)}\right),\\[6pt]
  F_{2} &=& (-1)^{r-1}\Delta^{r} \log
  \left(1+\frac{f_{2}(n)}{f_{1}(n)}\right),\\[6pt]
  F_{3} &=& (-1)^{r-1}\Delta^{r} \log \left(1+\frac{\sum_{k\geq 3}^{\infty}f_{k}(n)}{f_{1}(n)+f_{2}(n)}\right).
\end{eqnarray*}
Let
\begin{equation}\label{geq0loghrs2}
G_{r}=F_{1}+F_{2}+F_{3}.
\end{equation}
To estimate $(-1)^{r-1}\Delta^{r}\log p(n)$, we shall give upper bounds for $H_{r}$ and $G_{r}$.
We first consider  $G_{r}$.

\begin{theo}\label{estg}
For $n \geq 50$, we have 
\begin{equation}\label{gestmm1}
|G_{r}|<5\cdot 2^{r+\frac{1}{2}}e^{-\frac{\mu(n)}{2}}.
\end{equation}
\end{theo}

To prove Theorem \ref{estg}, we recall a monotone property
of the ratio of two power series, see Ponnusamy and Vuorinen \cite{pv1}. We also need a lower bound and an upper bound on the ratio of
 $L_{\nu}(x)$ and $L_{\nu}(y)$, which can be deduced from
known bounds  on the ratio of two  modified Bessel functions of the first kind.

\begin{prop}\label{logpvm}
Suppose that the power series
\begin{equation*}
f(x)=\sum_{m=0}^{\infty}\alpha_{m}x^{m} \quad \textrm{and}\quad g(x)=\sum_{m=0}^{\infty}\beta_{m}x^{m}
\end{equation*}
both converge for $|x|<\infty$ and $\beta_{m}>0$ for all $m>0$. Then the function $\frac{f(x)}{g(x)}$ is strictly decreasing for $x>0$ if the sequence $\{\alpha_{m}/\beta_{m}\}_{m=0}^{\infty}$ is strictly decreasing.
\end{prop}

Let $I_{\nu}(x)$ be the modified Bessel function of the first kind as given by
\begin{equation*}
  I_{\nu}(x)=\left(\frac{x}{2}\right)^{\nu}\sum_{m=0}^{\infty}
  \frac{\left(\frac{x^{2}}{4}\right)^{m}}{m!\Gamma(m+\nu+1)},
\end{equation*}
see Watson \cite{watson1}.
It is known that for $\nu \geq 1/2$ and $0<x<y$, $I_{\nu}(x)$ increases with $x$ and
\begin{equation*}
e^{x-y}\left(\frac{x}{y}\right)^{\nu}<\frac{I_{\nu}(x)}{I_{\nu}(y)}
<e^{x-y}\left(\frac{y}{x}\right)^{\nu},
\end{equation*}
see Baricz \cite[inequalities 2.2 and 2.4]{bondli}.
For $x>0$, from \eqref{deflnux1} we see that $L_{\nu}(x)$ can be expressed
by  $I_{\nu}(x)$,
\begin{equation*}
L_{\nu}(x)=x^{-\nu/2}I_{\nu}(2\sqrt{x}).
\end{equation*}
Thus the above properties of $I_{\nu}(x)$ can be restated in terms of $L_{\nu}(x)$.

\begin{prop}\label{logpropb1}
For $\nu \geq 1/2$ and $0<x<y$, we have
\begin{equation*} 
  e^{2\sqrt{x}-2\sqrt{y}}<\frac{L_{\nu}(x)}{L_{\nu}(y)}<e^{2\sqrt{x}
  -2\sqrt{y}}\left(\frac{y}{x}\right)^{\nu}.
\end{equation*}
\end{prop}

We are now ready to prove Theorem \ref{estg}.

\noindent {\emph{Proof of Theorem \ref{estg}}.} Since  $|G_{r}|\leq |F_{1}|+|F_{2}|+|F_{3}|$,
in order to estimate $G_r$, we shall estimate $|F_{1}|$, $|F_{2}|$ and $|F_{3}|$.
It follows from \eqref{amr2}  that
\begin{eqnarray*}
  |f_{k}(n)|&=&2\pi \left(\frac{\pi}{12}\right)^{3/2}|A_{k}(n)|k^{-5/2}
  L_{3/2}\left(\frac{\mu(n)^2}{4k^2}\right)\\[6pt]
  &\leq&4\pi \left(\frac{\pi}{12}\right)^{3/2}k^{-7/4}
      L_{3/2}\left(\frac{\mu(n)^2}{4k^2}\right),
\end{eqnarray*}
which yields that
\begin{equation}\label{loghrs4}
\sum_{k=3}^{\infty}|f_{k}(n)|\leq 4\pi \left(\frac{\pi}{12}\right)^{3/2}\zeta(7/4)L_{3/2}\left(\frac{\mu(n)^2}{4k^2}\right),
\end{equation}
where $\zeta(x)$ is the Riemann zeta function.
For convenience, we denote by $g(n)$ the right hand side of the above inequality,
so that \eqref{loghrs4} becomes
\begin{equation}\label{loghrs5}
  \sum_{k=3}^{\infty}|f_{k}(n)|\leq g(n).
\end{equation}
To estimate $F_{1}$, $F_{2}$ and $F_{3}$, we shall make use of the
 monotonicity of $\frac{\mu(n)+1}{\mu(n)-1}e^{-2\mu(n)}$, $\frac{|f_{2}(n)|}{f_{1}(n)}$ and $\frac{g(n)}{f_{1}(n)-|f_{2}(n)|}$.
It is easily seen that $\frac{\mu(n)+1}{\mu(n)-1}e^{-2\mu(n)}$ decreases with $n$ for $n\geq 1$, since $\frac{y+1}{y-1}e^{-2y}$ decreases with $y$ for $y>0$ and $\mu(n)$ increases with $n$.
By \eqref{askm1}, we have
\begin{equation*}
  \frac{|f_{2}(n)|}{f_{1}(n)}=
  \frac{L_{3/2}(\mu^{2}(n)/16)}{2^{5/2}L_{3/2}(\mu^{2}(n)/4)}.
\end{equation*}
The ratio of coefficients of $x^m$ in $L_{3/2}(\mu^{2}(n)/16)$ and $L_{3/2}(\mu^{2}(n)/4)$ is $\frac{4^{m}}{16^{m}}$. By Proposition \ref{logpvm}, we see that $\frac{L_{3/2}(y/16)}{L_{3/2}(y/4)}$ decreases with $y$ for $y>0$.
 Notice that $\mu^{2}(x)$ increases with $x$ for $x\geq 1$. So $\frac{L_{3/2}(\mu^{2}(x)/16)}{L_{3/2}(\mu^{2}(x)/4)}$ decreases with $x$ for $x\geq 1$.
This implies that $\frac{|f_{2}(n)|}{f_{1}(n)}$ decreases with $n$.

Next we prove the monotonicity of $\frac{g(n)}{f_{1}(n)-|f_{2}(n)|}$.
Recall that
\begin{equation*}
 \frac{g(n)}{f_{1}(n)-|f_{2}(n)|} =  \frac{2\zeta(7/4)L_{3/2}(\mu^{2}(n)/36)}{L_{3/2}
  (\mu^{2}(n)/4)-2^{-5/2}L_{3/2}(\mu^{2}(n)/16)}.
\end{equation*}
The ratio of coefficients of $x^m$ in $L_{3/2}(y/36)$ and $L_{3/2}(y/4)-2^{-5/2}L_{3/2}(y/16)$ equals \[\frac{\frac{1}{36^{m}}}{\frac{1}{4^{m}}-\frac{1}{2^{5/2}16^{m}}},\]
which decreases with $m$ for $m\geq 0$. By Proposition \ref{logpvm}, we deduce that for $y>0$,
\begin{equation*}
\frac{L_{3/2}(y/36)}{L_{3/2}(y/4)-2^{-5/2}L_{3/2}(y/16)}
\end{equation*}
decreases with $y$. Hence $\frac{g(n)}{f_{1}(n)-|f_{2}(n)|}$ decreases with $n$ for $n\geq 1$.

Using the above monotone properties, we proceed to derive
 upper bounds for $|F_{1}|$, $|F_{2}|$ and $|F_{3}|$.
It is known that for $0<x<1$,
\begin{equation}\label{loglogeq1}
\log(1-x)\geq \frac{-x}{1-x},
\end{equation}
\begin{equation}\label{loglogeq2}
|\log(1\pm x)|\leq -\log (1-x),
\end{equation}
see also DeSalvo and Pak \cite{desp}.

We first estimate $F_{1}$. Since
\[\Delta^{r} f(n)=\sum_{k=0}^{r}(-1)^{r-k}{r \choose k}f(n+k),\]
we have
\begin{equation*}
F_{1}=\sum_{k=0}^{r}(-1)^{k+1}{r \choose k}\log \left(1+\frac{\mu(n+k)+1}{\mu(n+k)-1}e^{-2\mu(n+k)}\right).
\end{equation*}
It follows that
\begin{equation}\label{f1d1a1}
|F_{1}|\leq \sum_{k=0}^{r}{r \choose k}\log \left(1+\frac{\mu(n+k)+1}{\mu(n+k)-1}e^{-2\mu(n+k)}\right).
\end{equation}
By the monotonicity of $\frac{\mu(n)+1}{\mu(n)-1}e^{-2\mu(n)}$, we see that for $n\geq 1$ and $0\leq k \leq r$,
\begin{equation}\label{f1d1a2}
\log \left(1+\frac{\mu(n+k)+1}{\mu(n+k)-1}e^{-2\mu(n+k)}\right)\leq \log \left(1+\frac{\mu(n)+1}{\mu(n)-1}e^{-2\mu(n)}\right).
\end{equation}
Applying \eqref{f1d1a2} to \eqref{f1d1a1}, we find that for $n\geq 1$,
\begin{equation*}  
|F_{1}|\leq 2^{r}\log \left(1+\frac{\mu(n)+1}{\mu(n)-1}e^{-2\mu(n)}\right).
\end{equation*}
Since $\log(1+x)\leq x$ for $x\geq 0$, we see that for $n\geq 1$,
\begin{equation}\label{usef1}
|F_{1}|\leq 2^{r}\frac{\mu(n)+1}{\mu(n)-1}e^{-2\mu(n)}.
\end{equation}
To estimate $F_{2}$, we begin with the following expression
\begin{equation}\label{f2f2geq1}
F_{2}= \sum_{k=0}^{r}(-1)^{k+1}{r \choose k}\log\left(1+\frac{f_{2}(n+k)}{f_{1}(n+k)}\right).
\end{equation}
It follows from \eqref{geq0}   that
\[0<1-\frac{|f_{2}(n)|}{f_{1}(n)}<1.\]
Using \eqref{loglogeq2}, we find that for $n\geq 1$,
\begin{equation}\label{f2f2geq2}
\left|\log\left(1+\frac{f_{2}(n+k)}{f_{1}(n+k)}\right)\right|\leq -\log\left(1-\frac{|f_{2}(n+k)|}{f_{1}(n+k)}\right).
\end{equation}
Combining (\ref{f2f2geq1}) and (\ref{f2f2geq2}), we obtain that for $n\geq 1$,
\begin{equation*}
|F_{2}|\leq -\sum_{k=0}^{r}{r \choose k}\log\left(1-\frac{|f_{2}(n+k)|}{f_{1}(n+k)}\right).
\end{equation*}
Using the monotonicity of $\frac{|f_{2}(n)|}{f_{1}(n)}$, we see  that for $n\geq 1$,
\begin{equation*}
|F_{2}|\leq -2^{r}\log\left(1-\frac{|f_{2}(n)|}{f_{1}(n)}\right).
\end{equation*}
Hence, by \eqref{loglogeq1}, we obtain that for $n\geq 1$,
\begin{equation}\label{usef2}
|F_{2}|\leq 2^{r}\frac{|f_2(n)|}{f_1(n)-|f_2(n)|}.
\end{equation}

To estimate $F_{3}$, we use the following expression
\begin{equation}\label{f33est1}
F_{3}= \sum_{k=0}^{r}(-1)^{k+1}{r \choose k}\log\left(1+\frac{\sum_{k\geq3}^{\infty}f_{k}(n+k)}{f_{1}(n+k)+f_{2}(n+k)}\right).
\end{equation}
By Proposition \ref{logpropb1}, we find that for $n\geq 1$
\begin{equation}\label{f1f4a1}
  2^{-\frac{5}{2}}e^{-\frac{\mu(n)}{2}}<\frac{|f_{2}(n)|}{f_{1}(n)}<\sqrt{2}e^{-\frac{\mu(n)}{2}},
\end{equation}
and
\begin{equation}\label{f1f4a2}
  2\zeta(7/4)e^{-\frac{2\mu(n)}{3}}<\frac{g(n)}{f_{1}(n)}<54\zeta(7/4)e^{-\frac{2\mu(n)}{3}}.
\end{equation}
Consequently, for $n\geq 1$,
\begin{equation}\label{f1f4a3}
\frac{|f_{2}(n)|}{f_{1}(n)}+\frac{g(n)}{f_{1}(n)}<
\sqrt{2}e^{-\frac{\mu(n)}{2}}+54\zeta(7/4)e^{-\frac{2\mu(n)}{3}}.
\end{equation}
For $n\geq 50$, it can be checked that
\begin{equation}\label{f1f4a4}
\sqrt{2}e^{-\frac{\mu(n)}{2}}+54\zeta(7/4)e^{-\frac{2\mu(n)}{3}}<1.
\end{equation}
Combining \eqref{f1f4a3} and \eqref{f1f4a4}, we obtain that for $n\geq 50$,
\begin{equation*}
\frac{|f_{2}(n)|}{f_{1}(n)}+\frac{g(n)}{f_{1}(n)}<1,
\end{equation*}
or equivalently,
\begin{equation}\label{f1f4a5}
f_1(n)-|f_2(n)|-g(n)>0.
\end{equation}
Combining \eqref{loghrs5} and \eqref{f1f4a5}, we see that for $n\geq 50$,
\begin{equation*}
f_1(n)-|f_2(n)|-|\sum_{k\geq3}^{\infty}f_k(n)|>0,
\end{equation*}
which can be rewritten as
\begin{equation*}
1\geq 1-\frac{\left|\sum_{k\geq 3}^{\infty}f_{k}(n)\right|}{f_{1}(n)-|f_{2}(n)|}>0.
\end{equation*}
Thus, we can use \eqref{loglogeq2} to deduce that for $n\geq 50$,
\begin{equation}\label{f1f4a7}
\left|\log \left(1+\frac{\sum_{k\geq 3}^{\infty}f_{k}(n)}{f_{1}(n)+f_{2}(n)}\right)\right| \leq -\log\left(1-\frac{|\sum_{k\geq3}^{\infty}f_k(n)|}{f_{1}(n)-|f_{2}(n)|}\right).
\end{equation}
Since $-\log(1-x)$ is increasing for $x>-1$,
according to  \eqref{loghrs5} and \eqref{f1f4a7}, we deduce that for $n\geq 50$,
\begin{equation}\label{f1f4a8}
 -\log\left(1-\frac{|\sum_{k\geq3}^{\infty}f_k(n)|}{f_{1}(n)-|f_{2}(n)|}\right)
 <-\log\left(1-\frac{g(n)}{f_{1}(n)-|f_{2}(n)|}\right).
\end{equation}
Combining \eqref{f1f4a7} and \eqref{f1f4a8}, we see that for $n\geq 50$,
\begin{equation}\label{f1f4a9}
\left|\log \left(1+\frac{\sum_{k\geq 3}^{\infty}f_{k}(n)}{f_{1}(n)+f_{2}(n)}\right)\right| \leq -\log\left(1-\frac{g(n)}{f_{1}(n)-|f_{2}(n)|}\right).
\end{equation}
It follows  from (\ref{f33est1}) and (\ref{f1f4a9}) that for $n\geq 50$,
\begin{equation*}
|F_{3}|\leq -\sum_{k=0}^{r}{r\choose k}\log\left(1-\frac{g(n+k)}{f_{1}(n+k)-|f_{2}(n+k)|}\right).
\end{equation*}
Based on the monotonicity of $\frac{g(n)}{f_{1}(n)-|f_{2}(n)|}$, we find that for $n \geq 50$,
\begin{equation*}
|F_{3}|\leq -2^{r}\log\left(1-\frac{g(n)}{f_{1}(n)-|f_{2}(n)|}\right).
\end{equation*}
Hence, by \eqref{loglogeq1}, we obtain that for $n\geq 50$,
\begin{equation}\label{usef3}
|F_{3}|\leq 2^{r}\frac{g(n)}{f_{1}(n)-|f_{2}(n)|-g(n)}.
\end{equation}
By Proposition \ref{logpropb1}, we see that for $n\geq 1$
\begin{equation}\label{gf1f4a3}
  2^{\frac{7}{2}}\zeta(7/4)e^{-\frac{\mu(n)}{6}}<\frac{g(n)}{|f_{2}(n)|}<
  27\sqrt{2}\zeta(7/4)e^{-\frac{\mu(n)}{6}}.
\end{equation}
In view of \eqref{usef1} and \eqref{f1f4a1}, we obtain that for $n\geq 50$,
\begin{equation}\label{estg1f1}
\frac{|F_{1}|}{F_{4}} < 2^{\frac{5}{2}}\frac{\mu(n)+1}{\mu(n)-1}e^{-\frac{3}{2}\mu(n)},
\end{equation}
where $F_{4}$ is defined by
\begin{equation*}
  F_{4}=2^{r}\frac{|f_{2}(n)|}{f_{1}(n)}.
\end{equation*}
As a consequence of \eqref{usef2} and \eqref{f1f4a1}, it can be checked that for $n\geq 50$,
\begin{equation}\label{estg1f2}
\frac{|F_{2}|}{F_{4}} < \frac{1}{1-\sqrt{2}e^{-\frac{\mu(n)}{2}}}.
\end{equation}
Applying \eqref{f1f4a1}, \eqref{f1f4a2} and \eqref{gf1f4a3} to \eqref{usef3}, we obtain that for $n\geq 50$,
\begin{equation}\label{estg1f3}
\frac{|F_{3}|}{F_{4}} < \frac{27\sqrt{2}\zeta(7/4)}{e^{\frac{\mu(n)}{6}}-
\sqrt{2}e^{-\frac{\mu(n)}{3}}-54\zeta(7/4)e^{-\frac{\mu(n)}{2}}}.
\end{equation}
Combining \eqref{estg1f1}, \eqref{estg1f2} and \eqref{estg1f3}, we conclude that for $n\geq 50$,
\begin{equation}\label{estgron1}
  |F_{1}|+|F_{2}|+|F_{3}|<5 F_{4}.
\end{equation}
It follows from \eqref{f1f4a1}  that for $n\geq 1$,
\begin{equation}\label{estgron2}
F_{4}<2^{r+\frac{1}{2}}e^{-\frac{\mu(n)}{2}}.
\end{equation}
Thus (\ref{estgron1}) and (\ref{estgron2}) lead to an upper bound for $|F_{1}|+|F_{2}|+|F_{3}|$.
This completes the proof. \qed

To prove Theorem \ref{dprm1}, we still need to estimate $H_r$ and we shall use two
relations due to  Odlyzko \cite{amo} on the relations between the higher order differences and derivatives.

\begin{prop}\label{abc}
 Let $r$ be a positive integer. Suppose that $f(x)$ is a function with infinite continuous derivatives for $x\geq 1$, and $(-1)^{k-1}f^{(k)}(x)>0$ for $k\geq 1$. Then for $r>1$,
 \begin{equation*}
  (-1)^{r-1} f^{(r)}(x+r)\leq (-1)^{r-1}\Delta^{r}f(x) \leq (-1)^{r-1}f^{(r)}(x).
 \end{equation*}
\end{prop}

\noindent{\emph{Proof of Theorem \ref{dprm1}}.}
First, we treat the case $r=1$, which states that for $n\geq 12$,
\begin{equation} \label{th1h1}
 \Delta \log p(n)< \log \left(1+\frac{\sqrt{6}\pi}{6\left(n+1\right)^{1/2}}\right).
\end{equation}
Since we have estimated $|G_{r}|$, we only need to estimate $H_{r}$ for $r=1$. By Proposition \ref{abc}, we have
\begin{equation}\label{h1h1f1}
H_{1}\leq \frac{2\pi}{\sqrt{24n-1}}-\frac{36}{24(n+1)-1}+\frac{12}{(24n-1)(1-\frac{6}{\pi \sqrt{24n-1}})}.
\end{equation}
We claim that for $n\geq 50$,
\begin{equation}\label{gh1h1f1}
H_{1}< \frac{\sqrt{6}\pi}{6\left(n+1\right)^{1/2}}
-\frac{1}{n+1}+\frac{5}{4\left(n+1\right)^{3/2}}.
\end{equation}
We proceed to estimate each term of the right hand side of \eqref{h1h1f1}. For the first term, we need to show that for $n\geq 50$,
\begin{equation}\label{h1h1f2}
\frac{2\pi}{\sqrt{24n-1}}<\frac{\sqrt{6}\pi}{6\left(n+1\right)^{1/2}}-\frac{3}{2(n+1)}.
\end{equation}
Setting $x =\frac{25}{24(n+1)}$, $\alpha=1/2$ and $c=\frac{1}{48}$, for $n\geq 50$, we have $0<x<c<1$ and $\alpha\geq \frac{1}{2}$. It follows from (\ref{gongs1}) that for $n\geq 50$,
\begin{eqnarray*}
  \frac{2\pi}{\sqrt{24n-1}}&=&\frac{2\pi}{\sqrt{24}\left(n+1\right)^{1/2}
  \left(1-\frac{25}{24(n+1)}\right)^{1/2}} \nonumber \\[6pt]
  &\leq& \frac{2\pi}{\sqrt{24}\left(n+1\right)^{1/2}}\left(
  1+\left(\frac{48}{47}\right)^{3/2}\frac{25}{48(n+1)}\right).
\end{eqnarray*}
This proves \eqref{h1h1f2}.

For the second term of the right hand side of \eqref{h1h1f1},  for $n\geq 50$,
we have \begin{equation}\label{h1h1f3}
-\frac{36}{24(n+1)-1}<-\frac{3}{2(n+1)}.
\end{equation}
For the last term of the right hand side of \eqref{h1h1f1}, using the same argument as in the derivation of \eqref{22s1}, we obtain that for $n\geq 50$,
\begin{equation}\label{h1h1f4}
\frac{12}{(24n-1)(1-\frac{6}{\pi \sqrt{24n-1}})}<\frac{1}{2(n+1)}
+\frac{1}{2\left(n+1\right)^{3/2}}.
\end{equation}
Combining \eqref{h1h1f2}, \eqref{h1h1f3} and \eqref{h1h1f4}, we obtain \eqref{gh1h1f1}.

By the estimate of $H_1$ in \eqref{gh1h1f1} and the
estimate of $G_1$ in \eqref{gestmm1}, we obtain that for $n\geq 50$,
\begin{equation*}
\Delta \log p(n)<\frac{\sqrt{6}\pi}{6\left(n+1\right)^{1/2}}
-\frac{1}{n+1}+\frac{5}{4\left(n+1\right)^{3/2}}
+10\sqrt{2}e^{-\frac{\pi}{12}\sqrt{(24n-1)}}.
\end{equation*}
Notice that for $n\geq 200$,
\begin{equation*}
\frac{5}{4\left(n+1\right)^{3/2}}<\frac{12-\pi^{2}}{24(n+1)},
\end{equation*}
and for $n\geq 50$,
\begin{equation*} 10\sqrt{2}e^{-\frac{\pi}{12}\sqrt{(24n-1)}}<\frac{12-\pi^{2}}{24(n+1)}.
\end{equation*}
Hence, for $n\geq 200$,
\begin{equation}\label{ggh1ga1}
\Delta \log p(n)<\frac{\sqrt{6}\pi}{6\left(n+1\right)^{1/2}}
-\frac{\pi^{2}}{12(n+1)}.
\end{equation}
Moreover, it can be easily checked that for $x>0$,
\[
x\left(1-\frac{x}{2}\right)< \log(1+x).
\]
Thus, for $n\geq 1$,
\[
\frac{\sqrt{6}\pi}{6\left(n+1\right)^{1/2}}-\frac{\pi^{2}}{12(n+1)}
<\log\left(1+\frac{\sqrt{6}\pi}{6\left(n+1\right)^{1/2}}\right).
\]
Combining the above relation and (\ref{ggh1ga1}), we reach \eqref{th1h1} for $n\geq 200$.

It can be checked that \eqref{th1h1} is valid for $12 \leq n \leq 200$,
and so Theorem \ref{dprm1} holds for $r=1$.

We now turn to the case $r\geq2$. We proceed to show that  there exists an integer $n(r)$
such that for $n\geq n(r)$,
\begin{equation}\label{drum0}
(-1)^{r-1}\Delta ^{r} \log p(n)< U_{r},
\end{equation}
where
\begin{equation*}
  U_{r}= \frac{\sqrt{6}\pi}{6} \left(\frac{1}{2}\right)_{r-1}\frac{1}{(n+1)^{r-\frac{1}{2}}}\left(1-\frac{\sqrt{6}\pi}{6} \left(\frac{1}{2}\right)_{r-1}\frac{1}{(n+1)^{r-\frac{1}{2}}}\right).
\end{equation*}
Since $x(1-x)<\log(1+x)$ for $x>0$, we have that for $n\geq 1$,
\begin{equation*}
  U_{r}<\log \left(1+\frac{\sqrt{6}\pi}{6} \left(\frac{1}{2}\right)_{r-1}\frac{1}{(n+1)^{r-\frac{1}{2}}}\right).
\end{equation*}
Thus \eqref{drum0} implies Theorem \ref{dprm1} for $r\geq 2$.

By \eqref{loghrs2}, we see that for $n\geq 1$,
\begin{equation*}
(-1)^{r-1}\Delta^{r}\log p(n)\leq H_r+|G_r|.
\end{equation*}
To prove \eqref{drum0}, it suffices to show that for $n\geq n(r)$
\begin{equation}\label{drumm1}
H_r+|G_r|<U_r.
\end{equation}
Since Theorem \ref{estg} gives the upper bound for $|G_r|$, we  need an upper bound for $H_r$. Recall that for $n\geq 1$,
\begin{equation}\label{husum1}
H_{r}=(-1)^{r-1}\Delta^{r}\left(-3\log \mu(n)
  +\log (\mu(n)-1)+\mu(n) \right).
\end{equation}
For $x\geq1$, write
\[\log (\mu(x)-1)=\log \mu(x)-\sum_{k=1}^{\infty}\frac{1}{k\mu(x)^{k}}.\]
By  exchanging  the order of two summations with one being finite,
it can be seen that  for $x\geq 1$,
\[
\Delta ^{r}\log (\mu(x)-1)=\Delta \log \mu(n)-\sum_{k}^{\infty}\Delta^{r}\left(\frac{1}{k\mu(n)^{k}}\right).
\]
Hence  \eqref{husum1} implies that for $n\geq 1$,
\begin{equation*}
H_r=(-1)^{r-1}\Delta^{r}\left(\mu(n)-2\log \mu(n)\right)-\sum_{k=1}^{\infty}(-1)^{r-1}\Delta^{r}
\left(\frac{1}{k\mu(n)^{k}}\right).
\end{equation*}
The $r$th derivatives of $\mu(x)=\frac{\pi}{6}\sqrt{24x-1}, \log \mu(x)$ and $\mu(x)^{-k}$ are given as follows,
\begin{eqnarray*}
\mu^{(r)}(x)&=&\frac{(-1)^{r-1}(\frac{1}{2})_{r-1}24^{r}\pi}{12(24x-1)^{r-\frac{1}{2}}}, \\[6pt]
\log^{(r)}(\mu(x))&=&\frac{(-1)^{r-1}(r-1)!24^{r}}{(24x-1)^{r}}, \\[6pt]
\left(\frac{1}{\mu^{k}}\right)^{(r)} &=& \left(\frac{k}{2}\right)_{r} \frac{(-144)^{r}}{\pi^{k}(24x-1)^{\frac{k}{2}+r}}.
\end{eqnarray*}
Therefore, the functions $\mu(x)=\frac{\pi}{6}\sqrt{24x-1}, \log \mu(x)$ and $-\mu(x)^{-k}$ satisfy the conditions of Proposition \ref{abc} for  $r\geq 1$ and $k\geq 1$. Hence,
\begin{eqnarray}\label{upbnd1}
   H_{r} &\leq& \frac{(\frac{1}{2})_{r-1}24^{r}\pi}{12(24n-1)^{r-\frac{1}{2}}}
  -\frac{(r-1)!24^{r}}{(24(n+r)-1)^{r}} \nonumber \\[6pt]
   &&\quad + \sum_{k=1}^{\infty}\left(\frac{k}{2}\right)_{r} \frac{144^{r}}{k\pi^{k}(24n-1)^{\frac{k}{2}+r}}.
\end{eqnarray}
Rewrite the right hand side of \eqref{upbnd1} as
\begin{align}\label{upbnd2}
&\frac{(\sqrt{6}\pi\frac{1}{2})_{r-1}}{(n+1)^{r-\frac{1}{2}}
\left(1-\frac{25}{24(n+1)}\right)^{r-\frac{1}{2}}} -\frac{(r-1)!}{(n+1)^{r}\left(1-\frac{24r-25}{24(n+1)}\right)^{r}} \notag \\[6pt]
&\quad +\sum_{k=1}^{\infty}\left(\frac{k}{2}\right)_{r}
\frac{6^{r}}{k\pi^{k}24^{\frac{k}{2}}\left(n+1\right)^{\frac{k}{2}+r}
\left(1-\frac{25}{24(n+1)}\right)^{\frac{k}{2}+r}}.
\end{align}
To bound the first term of \eqref{upbnd2}, we claim that for $n\geq 48r-3$,
\begin{equation}\label{ansupbnd1}
\frac{\sqrt{6}\pi(\frac{1}{2})_{r-1}}{6(n+1)^{r-\frac{1}{2}}
\left(1-\frac{25}{24(n+1)}\right)^{r-\frac{1}{2}}}\leq U_r
+\frac{ a_{1}}{(n+1)^{r+\frac{1}{2}}},
\end{equation}
where
\begin{equation*}
a_{1}= \left(\frac{1}{2}\right)_{r-1}\left(\frac{48}{47}\right)^{r+\frac{1}{2}}(2r-1)
\frac{25\pi}{24^{\frac{3}{2}}}
+\frac{\pi^{2}}{6}\left(\left(\frac{1}{2}\right)_{r-1}\right)^{2}
\frac{1}{(48r-2)^{r-\frac{3}{2}}}.
\end{equation*}
Setting $x =\frac{25}{24(n+1)}$, $\alpha=r-1/2$ and $c=\frac{1}{48}$,
for $n\geq 48r-3$, we have $0<x<c<1$ and $\alpha\geq \frac{1}{2}$.
 Invoking (\ref{gongs1}), we find that for $n\geq 48r-3$,
\begin{equation*}
  \frac{1}{\left(1-\frac{25}{24(n+1)}\right)^{r-1/2}}\leq  1+\left(\frac{48}{47}\right)^{r+1/2}\frac{25(2r-1)}{48(n+1)}.
\end{equation*}
This yields that for $n\geq 48r-3$,
\begin{align*}
&\frac{\sqrt{6}\pi(\frac{1}{2})_{r-1}}{6(n+1)^{r-\frac{1}{2}}
\left(1-\frac{25}{24(n+1)}\right)^{r-\frac{1}{2}}} \\[6pt]
&\quad \leq U_r+
\frac{\pi^{2}\left(\left(\frac{1}{2}\right)_{r-1}\right)^{2}}{6(n+1)^{2r-1}}+
\frac{25\pi(2r-1)\left(\frac{1}{2}\right)_{r-1}\left(\frac{48}{47}\right)^{r+\frac{1}{2}}
}{24^{3/2}\left(n+1\right)^{r+1/2}}.
\end{align*}
It is easily seen that for $n\geq 48r-3$,
\begin{equation*}
 \frac{1}{\left(n+1\right)^{2r-1}}\leq
 \frac{1}{\left(n+1\right)^{r+1/2}\left(48r-2\right)^{r-3/2}}.
\end{equation*}
So we arrive at \eqref{ansupbnd1}.

As for the second term of \eqref{upbnd2},  it can be shown that for $n\geq 48r-3$,
\begin{equation}\label{ansupbnd2}
-\frac{(r-1)!}{(n+1)^{r}\left(1-\frac{24r-25}{24(n+1)}\right)^{r}}\leq -\frac{(r-1)!}{(n+1)^{r}},
\end{equation}
or equivalently,
\[0<\frac{24r-25}{24(n+1)}< 1.\]
This can be easily justified since for $0\leq x<1$, $r\geq 2$ and $n\geq 48r-3$, \[\frac{1}{(1-x)^{r}}\geq 1 .\]

To estimate the last term of \eqref{upbnd2}, we aim to show that for $n\geq 48r-3$,
\begin{equation}\label{ansupbnd3}
\sum_{k=1}^{\infty}\left(\frac{k}{2}\right)_{r}
\frac{6^{r}}{k\pi^{k}24^{\frac{k}{2}}\left(n+1\right)^{\frac{k}{2}+r}
\left(1-\frac{25}{24(n+1)}\right)^{\frac{k}{2}+r}}\leq \frac{ a_{2}+a_{3}}{(n+1)^{r+\frac{1}{2}}},
\end{equation}
where
\begin{align*}
&a_{2}=\sum_{k=1}^{\infty}\left(\frac{k}{2}\right)_{r} \left(\frac{1}{48r-2}\right)^{\frac{k-1}{2}}
\frac{6^{k}}{k\pi^{k}24^{\frac{k}{2}}}, \\[6pt]
&a_{3}=\sum_{k=1}^{\infty}\left(\frac{k}{2}\right)_{r+1} \left(\frac{1}{48r-2}\right)^{\frac{k+1}{2}}
\left(\frac{48}{47}\right)^{\frac{k}{2}+r+1}\frac{25\cdot6^{k}(r+\frac{k}{2})}
{k\pi^{k}24^{\frac{k}{2}+1}}.
\end{align*}
Note that for given $r$, $a_2$ and $a_3$ are convergent. Setting $x =\frac{25}{24(n+1)}$, $\alpha=k/2+r$ and $c=\frac{1}{48}$, for $n\geq 48r-3$, we have $0<x<c<1$ and $\alpha\geq \frac{1}{2}$. By (\ref{gongs1}), we find that for $n\geq 48r-3$,
\begin{equation}\label{dpakhrg2}
  \frac{1}{\left(1-\frac{25}{24(n+1)}\right)^{r-1/2}}\leq  1+\left(\frac{48}{47}\right)^{k/2+r+1}\frac{25(2r+k)}{48(n+1)}.
\end{equation}
Clearly, for $n\geq 48r-3$ and $k\geq 1$,
\begin{eqnarray}
\frac{1}{\left(n+1\right)^{k/2+r}}&\leq& \frac{1}{\left(n+1\right)^{r+1/2}
\left(48r-2\right)^{\frac{k-1}{2}}}, \label{complexa1}\\
\frac{1}{\left(n+1\right)^{k/2+r+1}}&\leq&  \frac{1}{\left(n+1\right)^{r+1/2}
\left(48r-2\right)^{\frac{k+1}{2}}}. \label{complexa2}
\end{eqnarray}
Thus, \eqref{ansupbnd3} follows from \eqref{dpakhrg2}, \eqref{complexa1} and \eqref{complexa2}.

Combining \eqref{ansupbnd1}, \eqref{ansupbnd2} and \eqref{ansupbnd3}, we obtain that for $n\geq 48r-3$,
\begin{equation*}
  H_{r}(n)< U_r-\frac{(r-1)!}{(n+1)^{r}}+\frac{ a_{1}+a_{2}+a_{3}}{(n+1)^{r+\frac{1}{2}}}.
\end{equation*}
Let \[ u_{1}=\frac{4(a_{1}+a_{2}+a_{3})^{2}}{\left((r-1)!\right)^{2}}.\]
Notice that for given $r$ ,$a_1$ is a finite number and $a_2+a_3$ is convergent, so $a_1+a_2+a_3$ is a number for given $r$. It can be verified that for $n\geq u_1+1$,
\begin{equation*}
  \frac{ a_{1}+a_{2}+a_{3}}{(n+1)^{r+\frac{1}{2}}}<\frac{(r-1)!}{2(n+1)^{r}}.
\end{equation*}
Thus, for $n\geq \max\{48r-3, u_1+1\}$,
\begin{equation*}
H_{r}(n)<U_r-\frac{(r-1)!}{2(n+1)^{r}}.
\end{equation*}
Using the above inequality and \eqref{gestmm1}, we deduce that for $n\geq \max\{50, 48r-3, u_1+1\}$,
\begin{equation*}
H_{r}+|G_r|< U_{r}-\frac{(r-1)!}{2(n+1)^{r}}+5\cdot 2^{r+\frac{1}{2}}e^{-\frac{\mu(n)}{2}}.
\end{equation*}
Observe that for $n\geq 1$,
\begin{equation*}
\frac{1}{(n+1)^r}\geq \frac{\left(\frac{23}{48}\right)^{r}}{\left(n-\frac{1}{24}\right)^{r}}.
\end{equation*}
It follows that for $n\geq \max\{50, 48r-3, u_1+1\}$,
\begin{equation}\label{xyz1}
H_{r}+|G_r|< U_{r}-\frac{\left(\frac{23}{48}\right)^{r}(r-1)!}{2\left(n-\frac{1}{24}\right)^{r}}+5\cdot 2^{r+\frac{1}{2}}e^{-\frac{\mu(n)}{2}}.
\end{equation}
To deduce (\ref{drumm1}) from \eqref{xyz1}, we consider the following equation
\begin{equation}\label{expu1}
  \frac{\left(\frac{23}{48}\right)^{r}(r-1)!}{2\left(x-\frac{1}{24}\right)^{r}}
  =5\cdot 2^{r+\frac{1}{2}}e^{-\frac{\mu(x)}{2}}.
\end{equation}
Keep in mind that $\mu(x)$ is defined for $x\geq 1/24$. We claim that  equation \eqref{expu1} has two real roots. Recall that the Lambert $W$ function $W(z)$ is defined to be a function satisfying
\begin{equation}
W(z) e^{W(z)}=z,
\end{equation}
for any complex number $z$, see Corless, Gonnet, Hare, Jeffrey and Knuth \cite{lwfk}. So a solution of (\ref{expu1}) has the following form
\begin{equation*}
 x=\frac{1}{24}+\frac{6r^{2}}{\pi^2}
 \left(W\left(-\frac{\sqrt{46}\pi}{48r}
 \left(\frac{(r-1)!}{10\sqrt{2}}\right)^{\frac{1}{2r}}\right)\right)^{2}.
\end{equation*}
It is known that $W(z)$ is a multi-valued function. In particular, $W(z)$ has two real values $W_0(z)$ and $W_{-1}(z)$ for $-\frac{1}{e}<z<0$.
Using the following inequality, see Abramowitz and Stegun \cite{ahand},
\begin{equation}\label{estimatefac}
m!<\sqrt{2\pi}m^{m+\frac{1}{2}}e^{-m+\frac{1}{12m}},
\end{equation}
we see that for $r\geq 2$,
\begin{equation*}
\frac{\sqrt{46}\pi}{48r}
\left(\frac{(r-1)!}{10\sqrt{2}}\right)^{\frac{1}{2r}}<\frac{1}{e}.
\end{equation*}
Hence \eqref{expu1} has two real roots.
Let $u_{2}$ be the larger real root. Clearly, for sufficient large $x$,
\begin{equation*}
5\cdot 2^{r+\frac{1}{2}}e^{-\frac{\mu(x)}{2}}<
\frac{\left(\frac{23}{48}\right)^{r}(r-1)!}{2\left(x-\frac{1}{24}\right)^{r}}.
\end{equation*}
It follows that for $n\geq u_{2}+1$,
\begin{equation}\label{exepu1}
5\cdot 2^{r+\frac{1}{2}}e^{-\frac{\mu(n)}{2}}<\frac{\left(\frac{23}{48}\right)^{r}
(r-1)!}{2\left(n-\frac{1}{24}\right)^{r}}.
\end{equation}
Combining \eqref{xyz1} and \eqref{exepu1}, we conclude that
 \eqref{drumm1} holds  for $n\geq n(r)$, where \[ n(r)=\max\{50, 48r-3, u_{1}+1, u_{2}+1\}.\]
 This completes the proof for the case $r\geq 2$.
\qed

\section{The positivity of $(-1)^{r-1}\Delta^{r} \log p(n)$}

In this section, we prove the positivity of $(-1)^{r-1}\Delta^{r} \log p(n)$ for $r\geq 1$ and sufficiently large $n$. This is analogous to the positivity of  the
differences of the partition function conjectured by Good \cite{goodp} and proved by Gupta. The proof relies on the estimates of $H_r$ and $G_r$ in the previous section.

\begin{theo}\label{dprm2}
For each $r\geq 1$, there exists a positive integer $n(r)$ such that for $n\geq n(r)$,
\begin{equation}\label{thmmdp1}
(-1)^{r-1}\Delta ^{r} \log p(n)>0.
\end{equation}
\end{theo}

\proof The case $r=1$ is obvious since $p(n+1)>p(n)$ for $n\geq 1$.
For $r=2$,  DeSalvo and Pak \cite{desp} have shown that sequence $p(n)$ is log-concave for $n> 25$, or equivalently, for $n\geq 25$,
\begin{equation*}
-\Delta^{2}\log p(n)>0.
\end{equation*}
We now consider the case $r\geq 3$. Recall that
\begin{equation*}
(-1)^{r-1}\Delta ^{r}\log p(n)=H_r+G_r,
\end{equation*}
where $H_r$ and $G_r$ are given in (\ref{loghrs2})and (\ref{geq0loghrs2}).
Hence, we see that for $r\geq 1$,
\begin{equation}\label{dprme1}
  (-1)^{r-1}\Delta ^{r}\log p(n)\geq H_{r}-|G_{r}|.
\end{equation}
An   upper bound for $|G_r|$ has been given in Theorem \ref{estg}, so
we only need a lower bound for $H_r$.
By Proposition \ref{abc}, we find that
\begin{eqnarray}
   H_{r}&=&  (-1)^{r-1}\Delta ^{r}\left(\mu(n)-2\log \mu(n)-\sum_{k=1}^{\infty}\frac{1}{k\mu(n)^{k}}\right) \nonumber\\[6pt]
  &\geq& \frac{(\frac{1}{2})_{r-1}24^{r}\pi}{12(24(n+r)-1)^{r-\frac{1}{2}}}-\frac{(r-1)!24^{r}}{(24n-1)^{r}} \nonumber\\[6pt]
   &&\quad + \sum_{k=1}^{\infty}\left(\frac{k}{2}\right)_{r} \frac{ 144^{r}}{k\pi^{k}(24(n+r)-1)^{\frac{k}{2}+r}}. \label{geqff1}
\end{eqnarray}
The first term of the right hand side of (\ref{geqff1}) has the following lower
bound for $r\geq 48r-2$,
\begin{equation}\label{geq0fir1}
\frac{(\frac{1}{2})_{r-1}24^{r}\pi}{12(24(n+r)-1)^{r-\frac{1}{2}}}\geq
\frac{b_{1}}{n^{r-\frac{1}{2}}}-\frac{b_{2}}{n^{r}},
\end{equation}
where
\begin{eqnarray*}
b_{1}&=&\frac{\sqrt{6}\pi}{6} \left(\frac{1}{2}\right)_{r-1}, \\[6pt]
  b_{2}&=&\frac{\pi\sqrt{48r-2}}{24^{\frac{3}{2}}}
  \left(\frac{1}{2}\right)_{r}.
\end{eqnarray*}
Setting $x =\frac{24r-1}{24n}$ and $\alpha=r-1/2$, for $n\geq 48r-2$, we have $0<x<1$ and $\alpha\geq \frac{1}{2}$. It follows from (\ref{posalpha1}) that for $n\geq 48r-2$,
\begin{equation*}
\frac{1}{\left(1+\frac{24r-1}{24n}\right)^{r-\frac{1}{2}}}\geq
1-\frac{24r-1}{24n}\left(r-\frac{1}{2}\right),
\end{equation*}
or equivalently,
\begin{equation*}
\frac{(\frac{1}{2})_{r-1}24^{r}\pi}{12(24(n+r)-1)^{r-\frac{1}{2}}}\geq
\frac{\sqrt{6}\pi}{6} \left(\frac{1}{2}\right)_{r-1}\frac{1}{n^{r-\frac{1}{2}}}-
\frac{\sqrt{6}\pi}{6}\left(\frac{1}{2}\right)_{r}\frac{24r-1}{24n^{r+\frac{1}{2}}}.
\end{equation*}
Observing that for $n\geq 48r-2$,
\begin{equation*}
\frac{1}{n^{r+\frac{1}{2}}}\leq \frac{1}{\sqrt{48r-2}n^r},
\end{equation*}
we obtain \eqref{geq0fir1} for $n\geq 48r-2$.

For the second term of the right hand side of \eqref{geqff1}, we claim that for $n\geq 48r-2$,
\begin{equation}\label{geq0fir2}
\frac{(r-1)!24^{r}}{(24n-1)^{r}}\leq \frac{b_{3}}{n^{r}},
\end{equation}
where
\begin{equation*}
b_{3}= (r-1)! \left(1+\frac{r}{24}\left(\frac{1}{48r-2}\right)
  \left(\frac{48}{47}\right)^{r+1}\right).
\end{equation*}
Setting $x =\frac{1}{24n}$, $\alpha=r$ and $c=\frac{1}{48}$, for $n\geq 48r-2$, we have $0<x<c<1$ and $\alpha\geq \frac{1}{2}$. By (\ref{gongs1}), we see that
for $n\geq 48r-2$,
\begin{equation*}
\frac{1}{\left(1-\frac{1}{24n}\right)^{r}}\leq
1+\left(\frac{48}{47}\right)^{r+1}\frac{r}{24n}.
\end{equation*}
So we obtain  \eqref{geq0fir2}   for $n\geq 48r-2$.

Since the last term of the right hand side of \eqref{geqff1} is positive,
combining \eqref{geq0fir1} and \eqref{geq0fir2}, we deduce that for $n\geq 48r-2$,
\begin{equation}\label{dprme3}
H_{r}\geq \frac{b_{1}}{n^{r-\frac{1}{2}}}-\frac{b_{2}+b_{3}}{n^{r}}.
\end{equation}
To derive a simpler expression for a lower bound of  $H_r$, let \[m_{1}=\frac{4(b_{2}+b_{3})^2}{b_{1}^2}.\] Then we have that for $n\geq m_{1}+1$,
\begin{equation*}
\frac{b_{2}+b_{3}}{n^{r}}< \frac{b_{1}}{2n^{r-\frac{1}{2}}}.
\end{equation*}
It follows that for $n\geq \max\{48r-2,m_1+1\}$,
\begin{equation}\label{dprme4}
  H_r(n)> \frac{b_{1}}{2n^{r-\frac{1}{2}}}.
\end{equation}
By \eqref {dprme1} and \eqref{dprme4}, we find that for $n\geq \max\{50, 48r-2, m_1 +1\}$,
\begin{equation}\label{dprme5}
(-1)^{r-1}\Delta ^{r}\log p(n)>\frac{b_{1}}{2n^{r-\frac{1}{2}}}
-5\cdot 2^{r+\frac{1}{2}}e^{-\frac{\mu(n)}{2}}.
\end{equation}
Notice that for $r\geq 1$ and $n\geq 1$,
\begin{equation*}
\frac{1}{n^{r-\frac{1}{2}}}\geq \frac{\left(\frac{23}{24}\right)^{r-\frac{1}{2}}}
{\left(n-\frac{1}{24}\right)^{r-\frac{1}{2}}}.
\end{equation*}
Thus, for $n\geq \max\{50, 48r-2, m_1 +1\}$,
\begin{equation}\label{geq0dprme5}
(-1)^{r-1}\Delta ^{r}\log p(n)>\left(\frac{23}{24}\right)^{r-\frac{1}{2}}\frac{b_{1}}{2n^{r-\frac{1}{2}}}
-5\cdot 2^{r+\frac{1}{2}}e^{-\frac{\mu(n)}{2}}.
\end{equation}
To prove that the right hand side of \eqref{geq0dprme5}
is positive for sufficiently large $n$,   consider the following equation
\begin{equation}\label{expx1}
  \left(\frac{23}{24}\right)^{r-\frac{1}{2}}\frac{b_{1}}{2x^{r-\frac{1}{2}}}=5\cdot 2^{r+\frac{1}{2}}e^{-\frac{\mu(x)}{2}}.
\end{equation}
The solution of \eqref{expx1} can be expressed in terms of the Lambert
$W$ function, namely,
\begin{equation}
x=\frac{1}{24}+\frac{6\left(2r-1\right)^{2}}{\pi^{2}}
W\left(-\frac{\sqrt{46}\pi}{24(2r-1)}\left(\frac{\pi\left(\frac{1}{2}\right)_{r-1}}
{20\sqrt{6}}\right)^{\frac{1}{2r-1}}\right)^{2}.
\end{equation}
For $r\geq 1$, we have $\left(\frac{1}{2}\right)_{r}<r!$.
Using the estimate of $r!$ as given by  \eqref{estimatefac}, we obtain that for $r\geq 3$,
\begin{equation*}
-\frac{1}{e}<-\frac{\sqrt{46}\pi}{24(2r-1)}\left(\frac{\pi\left(\frac{1}{2}\right)_{r-1}}
{20\sqrt{6}}\right)^{\frac{1}{2r-1}}<0.
\end{equation*}
Thus \eqref{expx1} has two real roots. Let $m_{2}$ be the larger real root of  equation \eqref{expx1}. Clearly,  for sufficiently large $x$,
\begin{equation}
  \left(\frac{23}{24}\right)^{r-\frac{1}{2}}\frac{b_{1}}{2x^{r-\frac{1}{2}}}-5\cdot 2^{r+\frac{1}{2}}e^{-\frac{\mu(x)}{2}}>0.
\end{equation}
It follows that for $n\geq m_{2}+1$,
\begin{equation}\label{expx2}
  \left(\frac{23}{24}\right)^{r-\frac{1}{2}}\frac{b_{1}}{2n^{r-\frac{1}{2}}}-5\cdot 2^{r+\frac{1}{2}}e^{-\frac{\mu(n)}{2}}>0.
\end{equation}
Let \[n(r)= \max\{50, 48r-2, m_{1}+1, m_{2}+1\}.\]
Combining \eqref{dprme5} and \eqref{expx2}, we conclude that for $n\geq n(r)$,
\begin{equation}
  (-1)^{r-1}\Delta^{r}\log p(n)>0.
\end{equation}
This completes the proof.
\qed

\vspace{.3cm}

\noindent{\bf Acknowledgments.} This work was supported by the 973 Project, the PCSIRT Project of
the Ministry of Education and the National Science Foundation of China.

\end{document}